\documentclass[12pt, a4paper]{article}
\usepackage{amssymb,amsmath}
\usepackage[pdftex]{color}
\usepackage{pdfsync}
\newtheorem{theorem}{Theorem}
\newtheorem{proposition}{Proposition}

\newtheorem{corollary}{Corollary}
\newtheorem{lemma}{Lemma}
\newtheorem{remark}{Remark}
\newcommand{\tx}[1]{\mbox{\;{#1}\;}}

\numberwithin{equation}{section}
\numberwithin{lemma}{section}
\numberwithin{theorem}{section}
\numberwithin{definition}{section}
\numberwithin{remark}{section}
\numberwithin{corollary}{section}
\numberwithin{proposition}{section}
\begin{document}
\date{}
\title{Sublinear elliptic problems with a Hardy potential. }
\author{Catherine Bandle\thanks{Mathematisches Institut, Universit\"at
Basel, Rheinsprung 21, CH-4051Basel, Switzerland,  \texttt{catherine.bandle@unibas.ch}}\:
and
Maria Assunta Pozio\thanks{Dipartimento di Matematica,  Sapienza Universit\`{a} di Roma,
P.le  A. Moro 5, I-00185 Roma, Italy, \texttt{pozio@mat.uniroma1.it}}}
\pagestyle{myheadings}
\markright{\sc }
\maketitle

\abstract{In this paper we study the positive solutions of sub linear elliptic equations with a Hardy potential which is singular at the boundary. By means of ODE techniques a fairly complete picture of the class of radial solutions is given.
Local solutions with a prescribed growth at the boundary are constructed by means of contraction operators. Some of those radial solutions are then used to construct ordered upper and lower solutions in general domains. By standard iteration arguments the existence of positive solutions is proved. An important tool is the Hardy constant.}
\vspace{2mm}

\noindent
{\bf AMS Subject Classification}: 35J75, 
35B09, 
35B51, 
34B16. 
\smallskip

\noindent{\bf Key words}: Elliptic problems, Hardy potential, sub--linear forcing term, dead core solutions, boundary behavior.

\section{Introduction}
In this paper we study positive solutions of problems of the form
\begin{align}\label{original}
\Delta u +  \frac{\mu}{\delta(x)^2}u=u^p \tx{in} \Omega,
\end{align}
where $\mu \in \mathbb{R} \setminus \{0\}$, $\delta(x)$ is the distance of a point $x\in \Omega$ to the boundary,  $0<p<1$ and $\Omega \subset \mathbb{R}^N\,, \,N\ge 1$, is a bounded, smooth  domain. The expression
$$
\frac{\mu}{\delta(x)^2}=:V_\mu(x)
$$
is called the {\sl Hardy potential}.  In this type of problems there are two competing mechanisms, namely the nonlinear problem
\begin{align}\label{nonlinear}
\Delta u=u^p \tx{in} \Omega,
\end{align}
and the linear problem
\begin{align}\label{linear}
\Delta h+V_\mu(x)h=0 \tx{in} \Omega,
\end{align}
The problem\eqref{nonlinear} is well-understood cf. \cite{FrPh84}, \cite{BaSt84}. For any continuous function $\phi\geq 0$ it has a unique solution with $u=\phi$ on the boundary. Moreover if $\phi$  is small or if the domain is large the solutions have a {\sl dead core, i.e.} an open set $\omega \in \Omega$ where the solution vanishes identically. For the linear problem \eqref{linear} boundary values cannot be prescribed arbitrarily because of the singularity of the Hardy potential.

\medskip

\noindent The case  $p>1$ has been studied in \cite{BaMoRe08}. There among others, a partial classification of the solutions has been given.
It is based on the simple observation that the solutions of \eqref{original} are lower solutions for the linear problem,
and on some results of their local behavior near the boundary \cite{Ag83}. Another related study where the nonlinearity is the exponential function $e^u$ has been carried out in \cite{BaMoRe09}.
\medskip

\noindent It turns out that \eqref{original} has many solutions. We start with the investigation of radial solutions and provide a fairly complete picture of their structure. There are solutions whose boundary behavior is determined by the nonlinearity \eqref{nonlinear} and others by the Hardy potential \eqref{linear}.
In this latter case solutions have the same behavior as the harmonics of the $1$-- dimensional problem $h'' +V_{\mu}(x)h=0$ in $(-L,L)$.
Indeed the corresponding indicial equation is
\begin{equation}\label{beta}
\beta(\beta-1)+\mu=0.
\end{equation}
Hence positive harmonics near $x=-L$ and $x=L$ exist if and only if $\mu\leq 1/4$. Set
\begin{equation}\label{betapm}
\beta_\pm = \frac{1}{2}\pm \sqrt{\frac{1}{4}-\mu}.
\end{equation}
In this case $\delta(x) = L-|x|$ and the harmonics
are of the form ($ x \in (-L,0)$ or $ x \in (0,L)$ for any given constants $c_1,\,c_2 \in \mathbb{R}$)
$$
h(x) =
c_1 \delta^{\beta_+} +c_2 \delta^{\beta_-}\,,
$$
provided $\beta_-\neq \beta_+$. Otherwise
$$
h(x)=
c_1 \delta^{1/2} +c_2 \delta^{1/2} \log\frac{1}{\delta}.
$$
Observe that the derivative $h'$ is not in $L^2$ if $c_2$ is different from zero.
\smallskip

We will show that  there are only three possible boundary behaviors for the positive radial solutions, namely
\begin{align*}
1. \quad   &\lim_{\delta\to 0} \frac{u(\delta)}{\delta^{2/(1-p)}}=c',\tx{(nonlinear regime)}\\
2. \quad   &\lim_{\delta\to 0} \frac{u(\delta)}{\delta^{\beta_-}}=c_1,\tx{(linear singular regime)}\\
3. \quad   &\lim_{\delta\to 0} \frac{u(\delta)}{\delta^{\beta_+}}=c_2\tx{(linear regular regime)}.
\end{align*}
If the order of the linear regular regime is higher then the order of the nonlinear regime, only the second case occurs. We shall also prove the existence of local solutions with the boundary behavior described above.
These solutions are then used to construct upper and lower solutions in general domains.

\noindent An important tool for proving the existence of global solutions is the {\sl Hardy constant}. It is defined as
\begin{align}\label{Hardy}
C_H(\Omega)= \inf_{\phi\in W^{1,2}_0(\Omega)} \frac{\int_\Omega |\nabla \phi|^2\:dx}{\int_\Omega\delta^{-2}(x)\phi^2\:dx}.
\end{align}
It is well-known that $0<C_H\leq 1/4$ and $C_H(\Omega)=1/4$ for convex domains and for annuli if $N>2$, see Marcus, Mizel and Pinchover \cite{MaMiPi98}. If $N=2$ they proved that $C_H\to 0$ if the outer radius tends to infinity.  They also showed that for thin parallel sets $C_H=1/4$ and that the Hardy constant is attained if and only if $C_H<1/4$.

In our investigations  the following {\sl comparison principle} will play an important role:
\smallskip

{\sl  Let $\mu< C_H(\Omega)$ and $\omega\subseteq \Omega$. If $\Delta u+V_\mu u\geq 0$ in $\omega$ and $u\in W^{1,2}_0(\omega)$
 then $u\leq 0$ in $\omega$.}
 \medskip

 \noindent In fact  $u^+$ is an admissible function for \eqref{Hardy}.  Testing the inequality $\Delta u +V_\mu u \geq 0$ with $u^+$ we obtain $-\int_\Omega |\nabla u^+|^2\:dx + \mu\int_\Omega \frac{(u{^+)}^2}{\delta^2}\:dx \geq 0$. Hence $\mu \geq C_H(\Omega)$ which contradicts our assumption.

Our paper is organized as follows. We first study the radial solutions in balls and annuli. By means of ODE techniques we discuss the existence of local solutions with and without dead core and we show how to continue them globally. We then determine their asymptotic behavior near the boundary. At the end we prove the existence of positive solutions in arbitrary domains.


\section{Radial solutions, local behavior}
\subsection{Local solutions}\label{locsol}
In this section we study the radial solutions $u(r)$, $r=|x|$, of \eqref{original} in balls $B_R$ of radius $R$, centered at the origin, and in annuli $\mathcal{A}(r_0,R)=\{x: r_0<|x|<R\},\, r_0>0.$  They satisfy the ordinary differential equation  \begin{equation}\label{Er}
u'' +\frac{(N-1)}{r} u'
+\frac{\mu}{\delta(r)^2}u = u^p \tx{where} r
\in (0,R)  \tx{or}  r \in (r_0,R).
\end{equation}
Here $u'(r) :=\frac{d}{dr}u(r)$. It is well-known that problem \eqref{Er} with the initial conditions
\begin{align}
u(0)=u_0> 0,\: u'(0)=0 \label{ic1}\\
\nonumber \tx{or}\\
u(R_0)=u_0> 0,\: u'(R_0)=u_1\in \mathbb{R} \tx{for} R>R_0> r_0> 0\label{ic2}
\end{align}
has a unique local solution which is positive in a neighborhood of $r=0$ or of $R_0$, respectively.
Since the nonlinearity is not Lipschitz continuous at $u=0$,  the trivial solution is not the only solution
with $u(R_0)=0$ and $u'(R_0)=0$. In fact we shall prove that  there exists a local solution such that for a given $R_0\ge 0$ we have $u(R_0)=0$, $u'(R_0)=0$ and $u>0$ for $r>R_0$ and/or for $r<R_0$.
\subsubsection{Solutions with a dead core}
In our investigations there is a critical value of $\mu$ which will play an essential role. Define
\begin{align}\label{mucritical}
\mu^*:= \frac{2(p+1)}{(1-p)^2}.
\end{align}

\begin{lemma}\label{deadcore} (i) Let $R_0$ be a given point in $(r_0,R)$ in the case of an annulus, or in  $(0,R)$ in the case of a ball. Then in a small neighborhood of $R_0$ there exists a positive solution  of \eqref{Er} which is of the form $u(r)=|r-R_0|^{\frac{2}{1-p}}(c_p+w(r-R_0))$ and has the property that $u(R_0)=u'(R_0)=0$. Moreover $w(0)=0$ and $c_p=(\mu^*)^{\frac1{1-p}}$.

(ii) If  $R_0=r_0>0$ or $R_0=R$, then the same statement holds true provided $\mu>-\mu^*$. In this case  $c_p$ has to be replaced by  $c'=\big(\mu^*+\mu \big)^{\frac1{p-1}}$.

(iii) In the ball, near the origin, there exists a local solution of the form $u(r)= r^{\frac{2}{1-p}}(c'' +w(r))$ with $c''=\big(\mu^*+\frac{2(N-1)}{1-p} \big)^{\frac1{p-1}}$ and $w(0)=0$.
\end{lemma}
{\it Proof.}  Let us introduce in \eqref{Er} the new variable $d=r-R_0$. Then \eqref{Er} assumes the form
\begin{align*}
u''+\frac{N-1}{R_0+d}u' +\frac{\mu}{\delta^2} u =u^p \tx{in} (r_0-R_0, R-R_0).
\end{align*}
For simplicity we shall write $u(d)$ for $u(R_0+d)$.
Assuming that $u(0)=u'(0)=0$ we obtain after integration
$$
u(d)=\int_0^d\sigma(s)(u^p-\frac{\mu}{\delta^2}u)\left (\int_s^d\frac{dt}{\sigma(t)}\right)\:ds,
$$
where
$$
\sigma (d)= (R_0+ d)^{N-1}.
$$
Then
\begin{align}\label{Er4}
u(d)=& \int_0^d K_N(s,d;R_0)[u^p-\frac{\mu}{\delta^2}u] \,ds
\end{align}
 where
 \begin{align*}
 K_1&=d-s &\tx{if} N=1,\\
 K_2&=(R_0 + s)\ln\big(\frac{R_0+d}{R_0+s}\big) &\tx{if} N=2,\\
K_N&=\frac{R_0 +s}{N-2} [1-\big(\frac{R_0 +s}{R_0 + d} \big)^{N-2}]&\tx{if} N>2.
\end{align*}
The distance expressed in the variable $d$ becomes
$$
\delta(d)=
\begin{cases}
R-R_0-d \tx{if} R_0+d>(R+r_0)/2,\\
R_0-r_0+d \tx{if} R_0+d<(R+r_0)/2 .
\end{cases}
$$
From the Taylor expansion we obtain
\begin{equation}\label{KTaylor}
K_N(s,d;R_0)=d-s + O((d-s)^2)\,.
\end{equation}
Observe that \eqref{Er4} is also defined for negative $d$.
Set
\begin{equation}\label{Er5}
u(d):= |d|^{\frac2{1-p}}(c_p + w(d)))\,.
\end{equation}
By  \eqref{Er4} we have $w(d)=(Tw)(d)$ where
\begin{align}\label{Er6}
(Tw)(d):=
 \frac1{|d|^{\frac{2}{1-p}}} \int_0^dK_N(s,d) [ |s|^{\frac{2p}{1-p}}(c_p+w)^p-\frac{\mu}{\delta^2}|s|^{\frac2{1-p}}(c_p+w) ] \,ds-c_p.
 \end{align}
$Tw$ is well defined. Indeed by \eqref{KTaylor} we obtain for the lowest order term in the integral
\begin{equation}\label{cp}
\frac1{|d|^{\frac{2}{1-p}}} \int_0^d(d-s) |s|^{\frac{2p}{1-p}}\: ds = c_p^{1-p}\,.
\end{equation}
Next we want to show that $Tw$ has a fixed point. Fix $ \alpha$ such that
\begin{equation}\label{alfa}
p  < \alpha < 1\,,
\end{equation}
and define
\begin{equation}\label{Md}
M:= c_p\big(1-(\frac p{\alpha})^{\frac1{1-p}}\big) < c_p < 1\\
\nonumber \tx{and}
X:= \{w \in C^0([-d_0, d_0])\,: \,\,|w|_\infty \le M\}\,,
\end{equation}
where $d_0 \in (0, d_0^*]$ for some $d_0^*$ such that $\delta(\pm d_0^*)>0$. From the definition of $M$ it follows that
\begin{equation}\label{Malfa}
\alpha =\frac{p c_p^{1-p}}{(c_p-M)^{1-p}}<1\,.
\end{equation}
The following two properties hold:
\begin{itemize}
  \item[(i)] {\it $T$ is a contraction in $X$}.  A direct computation shows that for $d_0 \in (0, d_0^*]$
  $$
  0\le\int_0^dK_N(s,d)\:ds\leq C_0d^2 \tx{where  $C_0$ is independent of $d_0$}.
  $$
This together with \eqref{KTaylor}, \eqref{cp}, \eqref{Er6} and \eqref{Malfa} implies  that,  for sufficiently small $d_0$, and for a constant $C_1$ independent of $d_0$
\begin{equation}\label{Tcontr}
\begin{split}
|Tw_1 - Tw_2| \le  \frac{p c_p^{1-p}}{(c_p-M)^{1-p}}|w_1 - w_2|_\infty + C_1|d||w_1 - w_2|_\infty \\
= (\alpha + C_1 |d|) |w_1 - w_2|_\infty \le \frac{\alpha +1}2 |w_1 - w_2|_\infty\,.
\end{split}
\end{equation}

\item[(ii)] {\it $T:X\to X$}. \,\, For any $w \in X$ we have from the previous estimate
\begin{equation}\label{TXX}
\begin{split}
|Tw(\xi)| \le |Tw(\xi)-T0| +|T0|\\
\le \frac{\alpha + 1}2 |w|_\infty + C_1\,c_p |d| \\
\le \frac{\alpha + 1}2 M + C_1\,c _p|d| \le M\,,
\end{split}
\end{equation}
for $|d|\leq d_0$ sufficiently small.
\end{itemize}
Notice that by the special choice of $c_p$,  the fixed point satisfies $w(0)=0$ and consequently $u$ is positive in a neighborhood of $R_0$ .
\medskip

If $R_0=R$ or $r_0$, then $\delta(s) = |s|$, thus the linear term in \eqref{Er6} is of the form $\mu|s|^{\frac{2p}{1-p}}(c'+w)$. In order to have $w(0)=0$ we have to choose $c'$ suitably. With this change the remainder of the proof is the same as before. Similarly in the ball we have to adjust the constant if $R_0=0$. The details will  be omitted. This concludes the proof of the lemma. \hfill $\square$
\bigskip

In the previous lemma we have constructed a solution which vanishes together with its derivative at one point $R_0$. This solution gives rise to other solutions.
\begin{corollary}
For any $r_0<R_0'\leq R_0<R$, \eqref{Er} has a solution which is positive in $(R_0'-\epsilon,R_0')\cup (R_0,R_0+\epsilon)$ for $\epsilon>0$ sufficiently small and which vanishes in $[R_0',R_0]$.
We say that it has a dead core in $[R_0',R_0]$. Moreover there exist solutions vanishing in $(r_0,R_0)$ or in $(R_0,R)$ and positive in $(R_0,R_0+\epsilon)$ or in $(R_0-\epsilon,R_0)$.
\end{corollary}
\begin{corollary}\label{Coru=0}
Assume $\mu\in (-\mu^*, \frac14),\, \mu \not= 0$. If $u$ is a  local solution satisfying $\lim_{\delta \to 0} \frac{u(\delta)}{\delta^{\frac{2}{1-p}}}=0$ then $u\equiv 0$ in some neighborhood of the boundary.
\end{corollary}
{\sl Proof.} By contradiction suppose that there exists such a solution $u$ which is positive in $(0,\delta_0]$  ($\delta_0 >0$).\par
We assume $\delta_0>0$ so small that $C_H(\mathcal{A}(R- 2\delta_0,R))= \frac14$, and $C_H(\mathcal{A}(r_0,r_0 +2\delta_0))= \frac14$.
Then the maximum principle holds also if we are working in a larger annulus, if we deal with functions which belong to $W^{1,2}_0(\mathcal{A}(R- 2\delta_0,R))$ or $W^{1,2}_0(\mathcal{A}(r_0,r_0 +2\delta_0))$.\par
 First assume $\mu \in (0, \frac14)$. Let $\tilde{u}$ be the solution constructed in Lemma \ref{deadcore}, (ii). Since $\lim_{\delta \to 0} \frac{\tilde{u}(\delta)}{\delta^{\frac{2}{1-p}}}=c' =\big(\mu^*+\mu \big)^{\frac1{p-1}}>0$, we have
\begin{equation}\label{uutilde}
 u(\delta)< \tilde{u}(\delta)\,,\,\,\forall \delta \in (0, \delta_0],
\end{equation}
for a possibly smaller $\delta_0 >0$. For any $\epsilon \in (0, \delta_0)$, consider the function $\tilde{u}_\epsilon (\delta)$ such that $\tilde{u}_\epsilon (\delta)=0$ in $[0, \epsilon]$, and $\tilde{u}_\epsilon (\delta)= \tilde{u}(\delta -\epsilon)$ in $(\epsilon, \delta_0- \epsilon]$. By \eqref{uutilde} there exists  for sufficiently small $\epsilon$, a number $\delta_1 \in (0,\delta_0]$ such that $ \tilde{u}_\epsilon(\delta)<u(\delta)$ in $(0, \delta_1)$ and $ \tilde{u}_\epsilon(\delta_1)=u(\delta_1)$.
Since we have assumed that $\mu>0$ and since $ \tilde{u}_\epsilon$ belongs to $C^1([0,\delta_0])$, it can easily be  seen that $\tilde{u}_\epsilon$ is an upper solution of \eqref{original}, both near the inner or outer boundary. Then
$u-\tilde{u}_\epsilon=0$ at $\delta =0$ and $\delta = \delta_1$, $u-\tilde{u}_\epsilon>0$ and  satisfies $\Delta (u-\tilde{u}_\epsilon)+  \frac{\mu}{\delta(x)^2}(u-\tilde{u}_\epsilon)\geq u^p - \tilde{u}_\epsilon^p \ge 0$ for $0<\delta<\delta_1$.  This contradicts the maximum principle.
 Consequently $u$ vanishes in a neighborhood of  zero.\par
If $\mu \in (-\mu^*, 0)$, let $0<\epsilon<\delta_0$ and let $C$ be a positive number. Consider the function $z_\epsilon(\delta)=0$ in $[0, \epsilon]$, and $z_\epsilon(\delta) =C(\delta-\epsilon)^{\frac{2}{1-p}}$, for $\delta \in(\epsilon, \delta_0]$. $ z_\epsilon$ is $C^1([0,\delta_0])$ and for $\delta \in(\epsilon, \delta_0]$ it satisfies
\begin{equation}\label{Ez}
\begin {split}\Delta z_\epsilon + \frac{\mu}{\delta^2}z_\epsilon -z_\epsilon^p <  \Delta z_\epsilon -z_\epsilon^p\\
= C\delta^{\frac{2p}{1-p}}\left[\frac{2(p+1)}{(1-p)^2} +\frac{2(N-1)}{r(\delta)(1-p)}(\delta-\epsilon)
 -C^{p-1}\right],
 \end{split}
\end{equation}
where $r(\delta) =R-\delta$ at the outer boundary, $r(\delta) = r_0 + \delta$ at the inner boundary.
There exists a small positive constant $C_0$ depending only on $\delta_0$ such that the expression in the brackets of \eqref{Ez} is negative for all $\delta \in (\epsilon, \delta_0]$, hence $z_\epsilon$ is an upper solution for $C = C_0$. Because of our assumption we have $u\leq \frac{C_0}2\delta^{\frac 2{1-p}}$ in $(0,\delta_0]$ (for a possibly smaller $\delta_0$). Next we determine $\epsilon$ such that $z_\epsilon(\delta_0) \ge \frac{C_0}2 \delta_0^{\frac 2{1-p}}$. Then there exists $\delta_1 \le \delta_0$ such that $
z_\epsilon(\delta)<u(\delta)$ in $(0,\delta_1)$ and $z_\epsilon(\delta_1)=u(\delta_1)$. This is impossible by the comparison principle, as in the case of positive $\mu$.
Consequently $u$ vanishes in a neighborhood of  zero. \hfill $\square$


\subsubsection{Continuation of local solutions} Consider a local solution $u$ of the initial value problem \eqref{Er},
 \eqref{ic1} or \eqref{Er}, \eqref{ic2} respectively.  This solution can be continued up to the boundary unless it vanishes or blows up at an interior point. Blowup can be excluded because the nonlinearity is sub linear.
\medskip

Consider first a ball $B_R$. Assume that $u$ satisfies \eqref{Er}, \eqref{ic1} and $\mu<\frac{1}{4}$ . Then by the comparison principle stated in the Introduction $u$ cannot vanish at an inner point. Hence it can be continued as a global solution up to the boundary. By the same argument we can show that a solution with a dead core can be continued as a positive solution up to the boundary.  The positive solution to the left $(r<R_0)$ can be continued up to the origin but it is singular at the origin.
\smallskip

Consider now the solution of \eqref{Er}, \eqref{ic2} in an annulus and let  $\mu<C_H(\mathcal{A}(r_0,R))$. The solution can be continued at both sides until it vanishes or it reaches the boundary. By the comparison principle it cannot vanish at both sides at an interior point. Hence at least at one side it reaches the boundary. Thus a solution with a dead core can be continued as a positive solution which does not vanish at an inner point.

\subsection{Asymptotic behavior at the boundary}
In this section we assume that there exists a positive solution up to the inner or outer boundary in an annulus, and we want to determine the asymptotic behavior of $u$ as $r\to R$ or $r\to r_0$. The results in a ball coincide with those at the outer boundary of an annulus. \medskip

 {\sc Throughout this section we shall assume that $\mu<\frac{1}{4}$ and $\mu\neq 0$}.\medskip

The case $\mu=\frac{1}{4}$ can be treated similarly, but requires some further arguments and will therefore be omitted.

For this purpose we choose the distance from the boundary $\delta$ instead of $r$ as the new variable and
we write $u=\delta^\beta v$ where $\beta=\beta_+$ or $\beta_-$ defined in \eqref{betapm}, i.e.
$$
\beta_\pm= \frac{1}{2}\pm \sqrt{\frac{1}{4}-\mu}.
$$
From \eqref{Er}, for $\delta \in (0, \frac{R+r_0}{2})$
we obtain
\begin{align}\label{Ev0}
 v''+ \left (2\frac{\beta}{\delta}-\frac{N-1}{R-\delta}\right )v'-\beta\frac{N-1}{(R-\delta)\delta}v
=  v^p\delta^{\beta(p-1)} \tx{if} \delta = R-r\,,\\
\nonumber  v''+ \left (2\frac{\beta}{\delta}+\frac{N-1}{r_0+\delta}\right )v'+\beta\frac{N-1}{(r_0+\delta)\delta}v
=  v^p\delta^{\beta(p-1)} \tx{if} \delta = r-r_0\,.
\end{align}
These equations can be written in the form
\begin{align}
(\sigma_- v')'= \sigma_-\left(v^p\delta^{\beta(p-1)}+\beta\frac{N-1}{(R-\delta)\delta}v\right), \tx{where}
 \sigma_-(\delta)=\delta^{2\beta}(R-\delta)^{N-1},\label{Ev1}\\
 (\sigma_+ v')'= \sigma_+\left(v^p\delta^{\beta(p-1)}-\beta\frac{N-1}{(r_0+\delta)\delta}v\right), \tx{where}
 \sigma_+(\delta)=\delta^{2\beta}(r_0+\delta)^{N-1} \label{Ev1b}
\end{align}
\begin{lemma}\label{asympt.v1}
Let $v$ be a solution of \eqref{Ev1} with $\beta =\beta_->0$ or of \eqref{Ev1b} with $\beta=\beta_- <0$. Then
$$
\lim_{\delta\to 0} v(\delta)=v(0)<\infty.
$$
\end{lemma}
{\it Proof.} \:From the differential equations \eqref{Ev1} and \eqref{Ev1b} it follows immediately that for our particular choice of $\beta$, $v$ has no local maximum. It is therefore monotone near zero, hence there exists $\lim_{\delta \to 0} v(\delta)=v(0)$. Next we want to show that $v(0)<\infty.$ Suppose on the contrary that $v(0)=\infty$. Integration of \eqref{Ev1} yields
\begin{align*}
v(\delta)-v(\delta_0) +\sigma_-(\delta_0) v'(\delta_0)\int_\delta^{\delta_0} \sigma_-^{-1}\:ds =\\
 \int_\delta^{\delta_0}\sigma_- (v^ps^{\beta(p-1)}+\beta\frac{N-1}{(R-s)s}v)\:ds\int_\delta^s\sigma_-^{-1}\:d\xi.
\end{align*}
For $s\leq \delta_0$ we have since $\beta <1/2$
$$
\int_\delta^s \sigma_-^{-1}\:d\xi \leq \frac{s^{1-2\beta }}{(1-2\beta)(R-\delta_0)^{N-1}}.
$$
Since by assumption $v(\delta)$ is mono tone increasing near the origin
$$
v(\delta) \leq v(\delta_0) + c_1v'(\delta_0)\delta_0 + v^pg(\delta_0) + \frac{\beta(N-1)R^{N-1}}{(R-\delta_0)^{N}(1-2\beta)}\delta_0v(\delta),
$$
where $c_1$ and $g(\delta_0)$ are independent of $\delta$. We now choose $\delta_0$ so small that
$$
v(\delta) \leq v(\delta_0) + c_1v'(\delta_0)\delta_0 + v^pg(\delta_0) +\epsilon v(\delta) \tx{for $\epsilon <1$}.
$$
From here we deduce that $v(0)<\infty$. The same argument applies to the second statement.
\hfill $\square$
\medskip

If $\beta=\beta_-$ is of opposite sign the statement remains true but a different argument is required.
\begin{lemma} \label{asympt.v2} Let $v$ be a solution of \eqref{Ev1} with $\beta=\beta_- <0$ or of \eqref{Ev1b} with $\beta=\beta_- >0$. Then
$$
\lim_{\delta\to 0} v(\delta)=v(0)<\infty.
$$
\end{lemma}
{\it Proof.} \eqref{Ev1} and \eqref{Ev1b} imply that
$$
(\sigma_\pm v')'\leq \sigma_\pm v^p\delta^{\beta(p-1)}.
$$
Hence
\begin{align*}
v(\delta)\leq v(\delta_0)-\sigma_\pm(\delta_0)v'(\delta_0)\int_\delta^{\delta_0}\sigma_\pm^{-1}(s)
\:ds\\
+\int_\delta^{\delta_0}v^p\xi^{\beta(p-1)}\sigma_\pm\:d\xi\int_\delta^\xi\sigma_\pm^{-1}(t)\:dt,
\end{align*}
where $\sigma_-(\delta) = \delta^{2\beta}(R-\delta)^{N-1}$ and $\sigma_+(\delta) = \delta^{2\beta}(r_0+\delta)^{N-1}$.\par
\noindent
Since $1-2\beta>0$ and $\beta(p-1)+1>0$ it follows that
$$
v(\delta)\leq C_1 + C_2\int_\delta^{\delta_0}v^p\:d\xi,\tx{where $C_1, C_2$ are independent of $\delta$}.
$$
From this inequality we deduce that $v$ is uniformly bounded.
\bigskip

Next we want to show that $v(\delta)$ has a limit as $\delta$ tends to $0$.
\medskip

\noindent A. We first consider the equation \eqref{Ev1b} with $\beta >0$. After integration we obtain
\begin{align}\label{Eint1b}
\sigma_+(\delta) v'(\delta) - \sigma_+(\epsilon) v'(\epsilon) = \int_\epsilon^\delta \sigma_+\left(v^p s^{\beta(p-1)} - \beta\frac{(N-1)}{(r_0+s)s}v\right)\, ds
\end{align}
Notice that the right hand integral converges as $\epsilon \to 0$. We now  distinguish between two cases.
\medskip

1.\, $\lim_{\epsilon \to 0}\sigma_+(\epsilon) v'(\epsilon) = 0$. Then
$$
\sigma_+(\delta) v'(\delta) = \int_0^\delta \sigma_+\left(v^p s^{\beta(p-1)} - \beta\frac{(N-1)}{(r_0+s)s}v\right)\, ds.
$$
Since $v$ is bounded the following estimate holds true
$$
|v'(\delta)| \le \frac M{\delta^{2\beta}}\int_0^\delta s^{2\beta}\left(s^{\beta(p-1)} + \frac1{s}\right)\, ds \le c_1 \delta^{1+\beta(p-1)} +c_2\,,
$$
where $c_1,\, c_2$ are independent of $\delta$. Since $\beta <\frac{1}{2}$, $|v'|$ is bounded and hence $\lim_{\delta \to 0}v(\delta) = v(0)$.
\medskip

2.\, $\lim_{\epsilon \to 0}\sigma_+(\epsilon) v'(\epsilon)=L \not=0$. Then $v'(\delta) \to \pm \infty$ as $\delta \to 0$, depending on the sign of $L$. Again since $\beta_-<\frac{1}{2}$, $v(\delta)$ has a finite limit as $\delta\to 0$.
 \medskip

 B. Consider \eqref{Ev1} with $\beta<0$. As before we integrate \eqref{Ev1} and find

 $$
\sigma_- (\delta)v'(\delta)-\sigma_-(\epsilon)v'(\epsilon)=\int_\epsilon^\delta\sigma_-\left (v^ps^{\beta(p-1)}+\beta\frac{N-1}{(R-s)s}v\right)\:ds.
$$
Dividing this expression by $\sigma_-(\epsilon)$ we obtain the estimate
$$
|v'(\epsilon)|\leq c_1(\delta) \epsilon^{-2\beta} +c_2 |v|^p_\infty \epsilon^{-\beta(1-p)+1}+ c_3|v|_\infty\,,
$$
where $c_1,\, c_2,\,c_3$ are independent of $\epsilon$. Consequently $|v'|$ is bounded and
the  $\lim_{\delta \to 0}v(\delta)$ exists. This completes the proof of the lemma. \hfill $\square$
\medskip

As a consequence we can determine more precisely the behavior of $v(\delta)$ near zero.

\begin{proposition}\label{remark1} Let $v$ be a solution of \eqref{Ev1} or \eqref{Ev1b}. Then

(i) If  $\beta_- >0$ then  $\lim_{\epsilon \to 0} \sigma_\pm(\epsilon)v'(\epsilon)=L$.\\
Moreover if $L=0$ then
\begin{align}\label{as1}
v'(0)= \frac{N-1}{2R}v(0) \tx{at the outer boundary $r=R$}
\end{align}
and
\begin{align}\label{as2}
v'(0)= -\frac{N-1}{2r_0}v(0)\tx{at the inner boundary $r=r_0$}.
\end{align}
 \medskip

 (ii) If $\beta_-<0$ then  \eqref{as1} and \eqref{as2} hold.
  \end{proposition}
{\it Proof.}  By Lemmas \ref{asympt.v1}, \ref{asympt.v2}, we have that $v$ is continuous up to the boundary,
hence it is bounded. Near the outer boundary the function $v$ satisfies
\begin{equation}\label{intv}
\begin{split}
\sigma_- (\delta)v'(\delta)-\sigma_-(\epsilon)v'(\epsilon)\hspace{90pt}\\
=\int_\epsilon^\delta\sigma_-\left (v^ps^{\beta(p-1)}+\beta\frac{N-1}{(R-s)s}v\right)\:ds,
\end{split}
\end{equation}
where $\sigma_-(\delta) = \delta^{2\beta}(R-\delta)^{N-1}$.
If $\beta_->0$  the limit as $\epsilon$ tends to zero exists and is bounded, as we have already remarked. Let $\lim_{\epsilon \to 0} \sigma_-(\epsilon)v'(\epsilon)=L$ . If $L=0$ then
$$
\sigma_-(\delta)v'(\delta)=\int_0^\delta\sigma_-\left (v^ps^{\beta(p-1)}+\beta\frac{N-1}{(R-s)s}v\right)\:ds.
$$
Dividing by $\sigma_-(\delta)$ and applying the rule of Bernoulli l'Hospital we obtain \eqref{as1}.

If $\beta_-<0$ and $v(0)>0$, the integral at the right-hand side of \eqref{intv} becomes infinite as $\epsilon \to 0$. If we divide by
$\sigma_-(\epsilon)$ and apply the rule of Bernoulli l'Hospital the assertion follows. If $v(0)=0$ the same proof works if the integral at the right-hand side of \eqref{intv} diverges for $\epsilon \to 0$. If $v(0)=0$ and we have no information on the behavior of the integral as $\epsilon \to 0$, from \eqref{intv} we get
\begin{align*}
\frac{\sigma_- (\delta)v'(\delta)- \int_\epsilon^\delta\sigma_- v^ps^{\beta(p-1)}\:ds}{\sigma_-(\epsilon)}\le v'(\epsilon)\\
 \le \frac{\sigma_- (\delta)v'(\delta)+ |\beta| \int_\epsilon^\delta\sigma_- \frac{N-1}{(R-s)s}v\:ds}{\sigma_-(\epsilon)}\,.
\end{align*}
Taking the limit as $\epsilon \to 0$, we get $v'(0)=0$. Indeed the upper and lower bound are quotients. The denominator diverges, hence the limit is zero if the numerator is bounded. If the numerator is unbounded the application of the rule of Bernoulli l'Hospital gives the result.   This completes the proof for the outer boundary.\par The arguments in the case  of the inner boundary are exactly the same. \hfill $\square$
\bigskip

The main result of this section is summarized in
\begin{corollary} \label{Corb-}Suppose that $u$ exists and is positive up to the boundary. Assume $\mu<\frac{1}{4}, \, \mu \not=0$, and let $\beta_-=\frac{1}{2}-\sqrt{\frac{1}{4}-\mu}$. Then
$$
\frac{u(\delta)}{\delta^{\beta_-}}\to v(0)\tx{as} \delta \to 0\,.
$$
\end{corollary}
\medskip

Next we want to know more precisely what happens if $
\frac{u(\delta)}{\delta^{\beta_-}}\to 0\,\tx{as} \delta \to 0
$.
\begin{lemma}\label{le:b+} Assume $\mu <\frac14$, $\mu\neq 0$ and $
\frac{u(\delta)}{\delta^{\beta_-}}\to 0\,\tx{as} \delta \to 0
$.  Then there exists a nonnegative  constant $c\geq 0$ such that
$$
u(\delta) \leq c\delta^{\beta_+}
$$
in a neighborhood of the boundary.
\end{lemma}
{\it Proof.}  Set $u=\delta^{\beta}w$ where $\beta=\beta_+$. Then condition $\frac{u(\delta)}{\delta^{\beta_-}}\to 0\,\tx{as} \delta \to 0$ is equivalent to
\begin{equation}\label{w0}
  w(\delta)=\delta^{1-2\beta}v(\delta)\,,\quad\text{where}\quad v(\delta) \to 0\,\,\,\text{as}\,\,\delta \to 0.
\end{equation}
The functions $w$ and $v$ are solutions of \eqref{Ev0} with $\beta=\beta_+$ or $\beta_-$, respectively.
Integrating \eqref{Ev1}, or\eqref{Ev1b},respectively  in $[\epsilon, \delta]$ for $v$ replaced by $w$, we obtain
\begin{equation}\label{Ew0a}
\sigma_-(\delta) w'(\delta)-\sigma_-(\epsilon)w'(\epsilon)=\int_\epsilon^\delta \sigma_-(w^ps^{\beta(p-1)}+\beta\frac{N-1}{(R-s)s}w)\:ds\,
\end{equation}
and
\begin{equation}\label{Ew0b}
\sigma_+(\delta) w'(\delta)-\sigma_+(\epsilon)w'(\epsilon)=\int_\epsilon^\delta \sigma_+(w^ps^{\beta(p-1)}-\beta\frac{N-1}{(r_0+s)s}w)\:ds\,
\end{equation}
Since $\sigma_\pm \sim \delta^{2\beta}$ near zero it follows from our assumption that \newline
$\sigma(\epsilon) w^p\epsilon^{\beta(p-1)}\sim v^p\epsilon^{\beta(1-p)+p}$ and $\sigma(\epsilon)\frac{w}{\epsilon} \sim v$. Consequently the limit as $\epsilon \to 0$ exists at the right-hand sides of \eqref{Ew0a} and \eqref{Ew0b}, and thus by \eqref{w0} the limit as $\epsilon \to 0$ is finite. Hence there exists
$\lim_{\epsilon \to 0}\sigma(\epsilon) w'(\epsilon)=M$. If $M\neq 0$ then $w(\delta)\sim \delta^{1-2\beta}$
which is impossible since \eqref{w0} holds. Thus M=0, i.e. from \eqref{Ew0a} we get
\begin{equation}\label{Ew0c}
\sigma_-(\delta) w'(\delta)=\int_0^\delta \sigma_-(w^ps^{\beta(p-1)}+\beta\frac{N-1}{(R-s)s}w)\:ds\,.
\end{equation}
Hence $w'>0$ which implies,since $w$ is nonnegative, that $w$ is bounded as $\epsilon \to 0$.  From \eqref{Ew0b} we get
\begin{equation}\label{Ew0d}
\sigma_+(\delta) w'(\delta)=\int_0^\delta \sigma_+(w^ps^{\beta(p-1)}-\beta\frac{N-1}{(r_0+s)s}w)\:ds\,.
\end{equation}
We now integrate $w'(s)$ in the interval $[\epsilon, \delta]$.

\noindent If $\beta<1$ it follows from \eqref{w0}
that the integral converges for $\epsilon \to 0$, hence $w(0)$ exists and is finite. The case $\beta=1$ is excluded by our assumption $\mu\neq 0$.
\par
\noindent If $\beta>1$, we insert \eqref{w0} in \eqref{Ew0d} and neglect the positive term in the integral. Then
\begin{equation*}
 w'(\delta)\ge -c\delta^{-2\beta}\int_0^\delta s^{2\beta-1+1-2\beta}v(s)\:ds = -c\delta^{-2\beta}\int_0^\delta v(s)\:ds\,,
\end{equation*}
where $v(\delta) \to 0$ as $\delta \to 0$.
Integration from $\delta$ to $\delta_0<<1$ yields
$$
w(\delta_0)-w(\delta) \geq \frac{c}{2\beta-2}(\delta_0^{2-2\beta} -\delta^{2-2\beta}).
$$
Thus since $\beta>1$, $w(\delta)\le \delta^{2-2\beta} c_1$ for some positive constant $c_1$. Iterating this procedure of estimating $w'$ from below and $w$ from above, in a finite number of steps we get that $w$ is bounded.
This completes the proof.
 \hfill $\square$
\begin{remark}
If  $-\mu^*< \mu < \frac14$ then $\beta_+ < \frac2{1-p}$ and by Lemma \ref{deadcore}  there is a local solution which behaves like $c'\delta^{2/(1-p)}$ for a suitable $c'>0$. This solution is smaller than $\delta^{\beta_+}$and therefore $w(0)=0$ or equivalently $u(\delta)/\delta^{\beta_+} \to 0$ as $\delta \to 0$.

\noindent However if $\mu<-\mu^*$ or equivalently $\beta_+> \frac{2}{1-p}$, then no solution behaving like $c'\delta^{2/(1-p)}$ can exist. In fact if such a solution exists then it satisfies the assumptions of Lemma \ref{le:b+} and consequently $u\leq c\delta^{\beta_+}$. Hence
$$
0<c'=\lim_{\delta\to 0} \frac{u}{\delta^{2/(1-p)}} \leq \lim_{\delta\to 0}c\delta^{\beta_+-2/(1-p)}=0.
$$
This is impossible.
\end{remark}

If $\lim_{\delta\to 0} \frac{u}{\delta^{\beta_+} }=w(0)\neq 0$ then the  arguments developed in Proposition \ref{remark1} for the function $v$ apply also to the solution $w=u/\delta^{\beta_+}$.
\begin{proposition}\label{proposition2} Assume $ -\mu^*< \mu <\frac14$ and $\frac{u(\delta)}{\delta^{\beta_-}} \to 0$ as $\delta \to 0$. Then $w(\delta):=\frac{u(\delta)}{\delta^{\beta_+}} $ has a limit $w(0)$ for $\delta \to 0$ and if $w(0)\neq 0$ we get
$$ w(\delta) = w(0) + O(\delta^\alpha)\,,\,\,\tx{where}\,\,\alpha = \min\{1, 2-\beta_+(1-p)\}>0\,.$$
 Moreover if $-\frac{p}{(1-p)^2} < \mu $, we have $\alpha =1$ and
$$
w'(0)=\frac{N-1}{2R}w(0)\tx{at the outer boundary $r=R$}\,,
$$
$$
w'(0)=-\frac{N-1}{2r_0}w(0)\tx{at the inner boundary $r=r_0$.}
$$
If $w(0) =0$ we have $w(\delta) = o(\delta^\alpha)$.
\end{proposition}
{\it Proof.} We indicate the proof for the outer boundary. The statement for the inner boundary is proved in exactly the same way.

By the arguments given in the proof of Lemma \ref{le:b+} the function $w$ satisfies \eqref{Ew0c}.
If $-\frac{p}{(1-p)^2} \le \mu $, i.e. $\alpha =1$ we divide \eqref{Ew0c} by $\sigma_-$ and  apply Bernoulli l'Hospital's rule to get $w'(0)$. For $\mu = -\frac{p}{(1-p)^2}$
 the derivative involves an additional term. For $ -\mu^*< \mu <-\frac{p}{(1-p)^2}$, that is $0 <\alpha= 2-\beta_+(1-p) <1$, we divide \eqref{Ew0c} by $\sigma(\delta)$ and integrate in $[0, \delta]$ to get $w(\delta) -w(0)$ (cf. \eqref{Ew}). We easily see that the second term in the integral is higher order then the first one which is of order $\alpha$. By a simple analysis of it we get the conclusion
$$
\lim_{\delta \to 0}\frac{w(\delta)-w(0)}{\delta^\alpha}=\frac{w^p(0)}{(2-\beta_+(1-p))(1+\beta_+(1+p))}.
$$
 \hfill $\square$

The same type of argument as in Lemma \ref{le:b+} shows that there are no solutions which lie strictly between $c_1\delta^{\beta_+}$ and $c'\delta^{2/(1-p)}$.
 \begin{proposition}\label{proposition linear} Assume $-\mu^*<\mu<1/4$ or equivalently $\beta_+<\frac{2}{1-p}$. Then no solution exists for which $\lim_{\delta\to 0} u(\delta)\delta^{\frac{2}{p-1}}=\infty$ and $\lim_{\delta \to 0}u(\delta)\delta^{-\beta_+}=0$.
 \end{proposition}
 {\sl Proof.} By Proposition \ref{proposition2} we have $u(\delta)\delta^{-\beta_+} =o(\delta^\alpha)$ where $\alpha = \min\{1, 2-\beta_+(1-p)\}>0$.  Moreover $\beta_+< \beta_+ +\alpha \le \beta_+ +2-\beta_+(1-p) =2 + p\beta <\frac{2}{1-p}$ by our assumption on $\beta_+$.  Let $\epsilon$ be such that
 $$
 \beta_+(\epsilon) = 1/2+\sqrt{1/4-\mu +\epsilon}= \beta_+ +\alpha <\frac{2}{1-p} .
 $$
 By $\lim_{\delta\to 0} u(\delta)\delta^{\frac{2}{p-1}}=\infty$, there exists $\delta_0>0$ such that $u(\delta) \ge  \delta^\frac2{1-p}/\epsilon$ in $(0,\delta_0]$ and $C_H(\mathcal{A}(R- 2\delta_0,R))= \frac14 = C_H(\mathcal{A}(r_0,r_0 +2\delta_0))$ (cf. Corollary \ref{Coru=0}). Then
 $$
 -\frac{\mu}{\delta^2} u\leq \Delta u = u(u^{p-1}-\frac{\mu}{\delta^2}) \leq u\frac{\epsilon-\mu}{\delta^2} \tx{for} \delta >\delta_0.
 $$
 Let $h$ be a "small" harmonic satisfying $\Delta h +\frac{\mu-\epsilon}{\delta^2}h=0$ in a small neighborhood of the boundary. It behaves for $\delta$ near zero like $\delta^{\beta_+(\epsilon)}$. Since $h$ is defined up to a multiplicative constant, we can always assume that $h(\delta_0)=u(\delta_0)$. Remark that $h$ and $u$ are in $W^{1,2}_0$ in an neighborhood of the boundary, then the comparison principle applies and yields $u\geq h$. By the choice $\beta_+(\epsilon)= \beta_+ + \alpha$ we have $u(\delta) \ge c \delta^{\beta_+ + \alpha}$ for some positive $c$, which contradicts the result $u(\delta)\delta^{-\beta_+} =o(\delta^\alpha)$. Then the conclusion follows.
 \hfill $\square$

\subsection{Existence of local solutions at the boundary}
In this section we construct local solutions at the boundary points using the results of the previous sections. If $\mu >--\mu^*$, we know that a unique solution exists which grows at the  boundary like
$\delta ^{\frac2{1-p}}$ (cf. Lemma \ref{deadcore}, (ii)), i.e.  the non linear term is leading. Under the same condition on $\mu$, we also prove that there exist solutions which grow at the boundary like the harmonics.  It turns out that $\mu=-\frac{2(p+1)}{(p-1)^2}$ is a critical value.
\begin{lemma}\label{localregular}
Suppose that  $\beta= \beta_+$ and $-\mu^*< \mu < \frac14, \, \mu\not=0$.
For any positive constant $w(0)$ there exist near $r=R$ or $r=r_0$  a unique local solution of the form $u(\delta) = \delta^{\beta_+}w(\delta)$
where $w(\cdot)$ is continuous in $[0, \delta_0]$ (for some $\delta_0 \in (0, \frac{R-r_0}2)$) and
$$ w(\delta) = w(0) + O(\delta^\alpha)\,,\,\,\tx{where}\,\,\alpha = \min\{1, 2-\beta_+(1-p)\}>0\,.$$
Moreover $w(\cdot)$ behaves as described in Proposition \ref{proposition2}.
\end{lemma}
{\it Proof.}  Let us first consider the case  $r=R$. We write $u= \delta^{\beta_+}w$ and observe that $w$ satisfies equation \eqref{Ew0a} for $\beta=\beta_+$. We shall study the initial value problem \eqref{Ew0a} with $w(0) = w_0$ and $\lim_{\epsilon \to 0}\sigma_-(\epsilon)w'(\epsilon)=0$.
It can be transformed into the integral equation
\begin{align}\label{Ew}
w(\delta) = w(0)+ \underbrace{\int_0^\delta \sigma \left( s^{\beta(p-1)} w^p + \beta \frac{N-1}{(R-s)s}w\right)\left(\int_s^\delta\sigma^{-1}\,d\xi\right)\,ds}_{A(w)}\,,
\end{align}
where $\sigma(s) = \sigma_-(s)= s^{2\beta} (R-s)^{N-1}$.
Because of our assumption on $\mu$, we have $\beta_+(p-1)+1>-1$. Hence the integral exists for finite $w$. Straightforward computation shows that $A(w)$ is a contraction for small $\delta$. Hence there exists a fixed point $w$. Its behavior at zero follows from Proposition \ref{proposition2}. The same argument applies to the inner boundary. \hfill $\square$
\medskip
\begin{remark}\label{mumagg} The hypothesis $\mu>--\mu^*$ is necessary for the existence of solutions of order $\delta^{\beta_+}$. Indeed
the opposite condition $\mu \le --\mu^*$ is equivalent to $\beta_+(1-p)\geq 2$. Thus \eqref{Er} cannot be satisfied since for $u(r) = w(\delta)\delta^{\beta_+}$, $w(0) >0$, its left hand side depends on $\delta$ with order  higher then $\beta_+ -2$ and its right hand side is of order $p\beta_+ \le \beta_+-2$.
\end{remark}
\medskip

\begin{proposition} Assume $\mu<-\frac{2(p+1}{(1-p)^2}$ and $\lim_{\delta \to 0}\frac{u(\delta)}{\delta^{\beta_-}}=0.$ Then $u(x)\equiv 0 $ in a neighborhood of the boundary.
\end{proposition}
{\it Proof.} By contradiction suppose that such a local solution is positive in $[0,\delta_0)$ ($\delta_0>0$).\par
By lemma \ref{le:b+} our hypothesis on $u$ implies
\begin{equation}\label{uasb+}
u(\delta) \le c \delta^{\beta_+}\,,\,\,\,\, \delta \in [0, \delta _0]\,.
\end{equation}
For any $\epsilon \in [0, \delta_0)$  we define $z_\epsilon$ as in the proof of Corollary \ref{Coru=0}
for a constant $C=C_0$ such that $z_\epsilon$ is an upper solution of \eqref{original}. The hypothesis on $\mu$ gives $\beta_+ > \frac2{1-p}$,
hence by \eqref{uasb+} for a possibly smaller $\delta_0$ we get
$$
u(\delta) \le z_0(\delta)\,,\,\,\,\, \delta \in [0, \delta _0]\,.
$$
For  $\epsilon \in (0, \delta_0)$ sufficiently small we have $z(\delta_0)> u(\delta_0)$. Then $\delta_1\in (0,\delta_0)$
exists such that $z_\epsilon(\delta)<u(\delta)$ in $(0,\delta_1)$ and $z_\epsilon(\delta_1)=u(\delta_1)$.
As in the proof of Corollary \ref{Coru=0}, this is impossible by the comparison principle. \hfill $\square$
\medskip

In Lemma \ref{localregular} we constructed solutions that vanish on the inner or outer boundary and belong to the space $W^{1,2}$ in a neighborhood of it. Here  we prove the existence of "singular" local solutions near the boundary.
\begin{lemma}\label{localsingular1}
Let $\beta = \beta_-\in (0,1/2)$. For given $v(0)>0$ and $C\in \mathbb{R}$ there
exists near $r=R$ or $r=r_0$ a unique local solution of the following form respectively
\begin{align}\label{betaploc}
u(r) = \delta^{\beta_-}v(0)(1+ \frac{N-1}{2R}\delta +o(\delta)) + C \delta^{\beta_+}\,,\,\,\tx{if}\,\, \delta =R-r\,, \\
u(r) = \delta^{\beta_-}v(0)(1 -\frac{N-1}{2r_0} \delta + o(\delta)) + C \delta^{\beta_+}\,,\,\,\tx{if}\,\, \delta = r - r_0\,.
\end{align}
\end{lemma}
{\it Proof.} We look for a solution of the form $u=\delta^{\beta_-}v$ where $v$ satisfies
\eqref{Ev0}. At the outer boundary it can be written in view of Proposition \ref{remark1} as an integral equation of the form
\begin{align*}
v(\delta)-v(0)-L\int_0^\delta \sigma_-^{-1}ds=\\
\int_0^\delta\sigma_-\left(v^ps^{\beta_-(p-1)}+\beta_-\frac{N-1}{(R-s)s}v\right)
\int_s^\delta \sigma_-^{-1}\:d\xi \:ds\,,
\end{align*}
where $\sigma_-(s) = s^{2\beta_-}(R-s)^{N-1}$ and $L= C R^{N-1}(1-2\beta)$.
We can write $ v(\delta)=v(0)+L\int_0^\delta \sigma_-^{-1}ds + \delta \eta(\delta)$ and use a standard fixed point theorem to prove the existence of  $\eta$.  Moreover since $1-2\beta_- =\beta_+$ it follows that  $\delta^{\beta_-}L\int_0^\delta \sigma_-^{-1}ds = C \delta^{\beta_+} + \delta^{\beta_-}o(\delta)$.
Likewise we establish a solution at the inner boundary. \hfill $\square$
\medskip

The remaining case $\beta_-<0$ which requires a more subtle argument, is covered in the next lemma.
\begin{lemma}\label{LemmaMorels}
Assume $\beta=\beta_-<0$. Let  $v(0)$ be an arbitrary positive constant. Then there exists near $r=R$ or $r=r_0$ a one parameter family of local solutions of the form
$$
u(r)=\delta^{\beta_-}v(0)( 1+ \frac{N-1}{2R}\delta +\eta(\delta,C)) + C \delta^{\beta_+}
$$
or
$$
u(r)= \delta^{\beta_-}v(0)(1-\frac{N-1}{2r_0} \delta + \eta(\delta,C)) + C \delta^{\beta_+},
$$
respectively, where
\begin{align*}
\eta(\delta,C)=
\begin{cases}
o(\delta^{1-2\beta_-}) & \tx{if} -\frac{1}{2}<\beta_-<0,\\
K \delta^2 |\log(\delta)| +o(\delta^2|\log(\delta)|) & \tx{if} \beta_- =-\frac{1}{2},\\
K \delta^2 + o(\delta^2) & \tx{if} \beta_-<-\frac{1}{2}.
\end{cases}
\end{align*}
where $K$ is a constant which depends on the data of the problem but not on $v(0)$ neither on the parameter $C>0$ while the higher order terms do depend on both $v(0)$ and $C>0$. If $\beta \le -\frac12$ the term $C \delta^{\beta_+}$ would be included in the higher order term but we wrote it explicitly since this is the parameter on which solutions depend.
\end{lemma}
{\it Proof.}  As before we carry out the proof only for the outer boundary. Equation \eqref{Ev1} can be written as
$$
\frac{(v'(R-\delta)^{N-1})'}{(R-\delta)^{N-1}} + \frac{2}{\delta} \beta v' -\beta \underbrace{\frac{N-1}{(R-\delta)\delta}}_{\frac{N-1}{R}[\frac{1}{\delta} +\frac{1}{R-\delta}]}v=v^p\delta^{\beta(p-1)}.
$$
If we integrate this expression we get
\begin{align}\label{v2}
v'(\delta)(R-\delta)^{N-1} &-v'(0)R^{N-1}+\beta\int_0^\delta (R-s)^{N-1} \{2\frac{v'}{s} -\frac{N-1}{R}\frac{v}{s}\}\:ds
\nonumber \\& =
\int_0^\delta(R-s)^{N-1}\{ v^p s^{\beta(p-1)}+ \beta \frac{N-1}{R(R-s)}v\}\:ds.
\end{align}
Notice that that the integral at the left is singular without  additional conditions on $v$ and $v'$. Set
 $$
 v(\delta)=v(0)\left(1+ \frac{N-1}{2R} \delta +\eta(\delta)\right), \quad v'(\delta) =v(0)\left(\frac{N-1}{2R} + \eta'\right).
 $$
The behavior of the function $\eta$ near the origin will be specified later. For the moment we assume that all the integrals which appear in the calculations below are well-defined. Substituting $v$ and $v'$ we  get for the different expressions in the equation \eqref{v2}
\begin{align}\label{eq:10}
v'(\delta)(R-\delta)^{N-1}-v'(0)R^{N-1}= v(0)\eta' (R-\delta)^{N-1}\\
\nonumber  -v(0)\frac{(N-1)^2}{2}R^{N-3}\delta -\eta_1,
\end{align}

\begin{align}\label{eq:11}
-\beta \int_0^\delta (R-s)^{N-1}\{2\frac{v'}{s} -\frac{N-1}{R}\frac{v}{s}\}\:ds=-2\beta v(0) \int_0^\delta (R-s)^{N-1}\frac{\eta'}{s}\:ds\\
\nonumber +\beta v(0)\frac{(N-1)^2}{2}R^{N-3}\delta+\beta v(0)\frac{N-1}{R} \int_0^\delta (R-s)^{N-1}\frac{\eta}{s}\:ds+\eta_2
\end{align}

\begin{align}\label{eq:12}
\int_0^\delta(R-s)^{N-1}\{ v^p s^{\beta(p-1)}+ \beta \frac{N-1}{R(R-s)}v\}\:ds=\int_0^\delta(R-s)^{N-1} v^p s^{\beta(p-1)}\:ds\\
\nonumber  +\beta v(0)(N-1)R^{N-3}\delta +v(0) \frac{\beta(N-1)}{R}\int_0^\delta (R-s)^{N-2}\eta\:ds +\eta_3.
\end{align}
Here the functions $\eta_i$,  $i=1,2,3$, are of order $O(\delta^2)$   and  are independent of $\eta$ and $\eta'$. In the sequel we shall use the following notation:
\begin{align*}
A:= v(0)R^{N-3}\left\{ \frac{(N-1)^2}{2} +\beta(N-1) +\frac{\beta(N-1)^2}{2}\right\},\\
y_1:= v(0)(R-\delta)^{N-1} \eta' \tx{and} y_2:=v(0)(R-\delta)^{N-1} \eta,\\
H(\delta,y_2):=\int_0^\delta (R-s)^{N-1} v^p s^{\beta(p-1)} \:ds +\beta \frac{N-1}{R} \int_0^\delta (y_2+\frac{y_2}{s})\:ds+\sum_1^3\eta_i\,,
\end{align*}
where $v(s)$ is replaced by $v(0)(1+\frac{N-1}{2R}s + \frac{y_2(s)}{v(0)(R-s)^{N-1}})$. For the next arguments it will be important to keep in mind that $H$ is independent of $y_1$.
From \eqref{eq:10}, \eqref{eq:11} and \eqref{eq:12} it follows that
\begin{align}\label{y1}
y_1(\delta)=-2\beta \int_0^\delta \frac{y_1}{s}\:ds+A\delta +H(\delta,y_2).
\end{align}
Moreover
\begin{align}\label{y2}
y_2(\delta)=(R-\delta)^{N-1}\int_0^\delta \frac{y_1(s)}{(R-s)^{N-1}}\:ds.
\end{align}
If we set $\phi(\delta):= \int_0^\delta \frac{y_1}{s}\:ds$ then \eqref{y1} can be written as
\begin{align}\label{Ephi}
\phi'(\delta)\delta = -2\beta \phi(\delta) +A\delta +H(\delta,y_2)\,,
\end{align}
and solved by the variation - of - constants formula, as we will do below.
Clearly
\begin{align}\label{Y}
y_1= -2\beta\phi +A\delta +H(\delta,y_2)=:T(y_2).
\end{align}
We now replace $y_1$  in \eqref{y2} by $T(y_2)$ and we obtain a fixed point equation for $y_2$, namely
$$
y_2(\delta)= (R-\delta)^{N-1}\int_0^\delta\frac{T(y_2)}{(R-s)^{N-1}}\:ds=:\Theta(y_2).
$$
\smallskip

\noindent Next we want to show that in a properly chosen set, $\Theta(w)$ is a contraction. For this purpose we distinguish between three cases.
\medskip

(i) $-\frac{1}{2}<\beta<0$.
\smallskip

\noindent
Consider the Banach space $X:=\{w \in C([0,\delta_0]): |w|\leq M\delta^{-2\beta +1}, \delta \in [0,\delta_0]\}$ where $M$ and $\delta_0 \leq\frac{R-r_0}{2}$ are positive constants which will be determined later, and $\|w\|:=\sup\{|\frac{w_2}{\delta^{-2\beta +1}}|,\: \delta \in (0,\delta_0]\}$.  By the variation-of-constants formula \eqref{Ephi} gives
$$
\phi =
c\delta^{-2\beta} +\frac{A}{2\beta +1}\delta+\delta^{-2\beta}\int_0^\delta s^{2\beta-1}H(s,w)\:ds
$$
where $c$ is an arbitrary parameter.

First we estimate $|T(w)-T(\tilde w)| = |2\beta(\tilde \phi-\phi) +H(\delta,w)-H(\delta, \tilde w)|$ for $w$ and $\tilde w$ in $X$.
\smallskip

\noindent
Set for short
$v(\delta)=(R-\delta)^{1-N} (\gamma +w)$ where $\gamma =(R-\delta)^{N-1}v(0)(1+\frac{N-1}{2R}\delta)$ and similarly $\tilde v(\delta)$ with $w$ replaced by $\tilde w$. For fixed $M$ we can take a sufficiently small $\delta_0$ such that  $\gamma+w\ge (R-\delta)^{N-1}v(0)$.\par
Indeed $\delta_0 \le \frac{R-r_0}2$ implies $\gamma - (R-\delta)^{N-1}v(0) \ge \frac{(R+r_0)^{N-1}}{2^NR}v(0)\delta = c_0 \delta \ge M\delta^2 \ge |w|$ if
\begin{align}\label{M}
\delta_0 \le \frac{c_0}M
\end{align}
Then the following inequality holds
\begin{align}\label{vp}
|v^p-\tilde v^p|\leq\frac{p}{(R-\delta)^{(N-1)}v(0)^{1-p}}|w-\tilde w|.
\end{align}
Then
\begin{align*}
\big |\int_0^\delta(R-s)^{N-1}s^{\beta(p-1)}(v^p-\tilde v^p)\:ds\big| \leq c_1\delta^{\beta(p-3)+2}\|w-\tilde w\|,\\
|H(\delta,w)-H(\delta,\tilde w)| \leq c_1\delta^{\beta(p-3)+2}\|w-\tilde w\| + c_2\delta^{-2\beta+2}\|w-\tilde w\| \\+ c_3\delta^{-2\beta +1}\|w-\tilde w\|
 \leq c_4\delta^{-2\beta +1}\|w-\tilde w\|,
\end{align*}
where $c_i>0\, (i \in \mathbb{N} )$ stand for constants independent of $\delta$ all along this proof. Furthermore
\begin{align*}
|\phi(\delta)-\tilde \phi(\delta)|=|\delta^{-2\beta}\int_0^\delta s^{2\beta-1}(H(s,\tilde w)-H(s,w))\:ds|\\
\leq c_4\delta^{-2\beta+1}\|w-\tilde w\|.
\end{align*}
Therefore
$$
|T(w)-T(\tilde w)| \leq ( -2\beta c_4+ c_4)\delta^{-2\beta +1}\|w-\tilde w\|= c_4\delta^{-2\beta +1}\|w-\tilde w\|.
$$
Hence\par
$$
|\Theta(w)-\Theta(\tilde w)| \leq \frac{c_5}{-2\beta +2}\delta^{-2\beta+2}\|w-\tilde w\| $$
and \par
$\|\Theta(w)-\Theta(\tilde w)\|\leq c_6\delta_0\|w-\tilde w\|$.\par
For given $M$, $\delta_0$ can be chosen possibly smaller so that $\Theta(w)$ is a contraction and \eqref{M} holds.  It remains to show that $\Theta:X\to X$. For $w\in X$ we have by \eqref{Y}
$$
|T(w)| \leq |2\beta c |\delta^{-2\beta} + \frac{|A|}{2\beta +1}\delta + c_1 M\delta^{-2\beta +1},
$$
where $c_1$ is independent of M and $\delta_0$. Consequently
\begin{align*}
|\Theta(w)| \leq  \frac{2|\beta c|}{-2\beta +1}\delta^{-2\beta +1} +\frac{|A|}{2(2\beta +1)}\delta^2 +c_2\delta^{-2\beta +2}, \\
\|\Theta(w)\| \leq  \frac{2|\beta c|}{-2\beta +1} +\frac{|A|}{2(2\beta +1)}\delta^{1+2\beta} +c_2\delta.
\end{align*}
We now fix $M>-2\beta c$ and choose $\delta_0$ sufficiently small such that $\|\Theta(w)\| \leq M$. Notice that by decreasing $\delta_0$ the inequality \eqref{M} is not violated. Then $\Theta(w)$ is a contraction in $X$ and the conclusion follows. Indeed $C= \frac{-2\beta}{(1-2\beta)R^{N-1}}c \ge 0$ follows from the representation formula of the solution, as a fixed point, and from $\beta_+ = \beta_- + (1-2\beta_-)$.
This completes the proof for $\beta\in (-1/2,0)$
\par
 If $c=0$ we can carry out the same proof in the space
$X:=\{w \in C([0,\delta_0]): |w|\leq M\delta^2, \delta \in [0,\delta_0]\}$ with the norm $\|w\|:=\sup\{|\frac{w_2}{\delta^2}|,\: \delta \in (0,\delta_0]\}$. Here $M$ will be a constant close to $\frac{|A|}{2\beta +1}$
which is the leading term in $\|\Theta(w)\|$ if $A\neq 0$. For $c=0$ the solution is $C^2$ up to the boundary.
\medskip

(ii) $\beta =-\frac{1}{2}$.
\smallskip

\noindent In this case we have
$$
\phi=c\delta + A\delta |\log \delta|  +\delta\int_0^\delta Hs^{-2}\:ds.
$$
Here $-\beta+1 =2$ and we argue exactly as before if $A= 0$, i.e. $N=1$ or $N=3$. Otherwise the logarithmic term prevails. We then take $|w|\leq M \delta^2\log(1/\delta)$ and $\|w\|:=\sup\{|\frac{w_2(\delta)}{\delta^2|\ln \delta|}|\:, \,\,\delta \in (0,\delta_0]\}$. It turns out that for small $\delta_0$, $\Theta(w)$ is a contraction which maps the ball $\{|w|\leq M\delta^2|\log\delta|\}$ into itself. It has therefore a fixed point. The details will be omitted.
\medskip

 (iii) $\beta<-\frac{1}{2}$.
\smallskip

The function $\phi$ defined before is in general not defined for $\delta=0$ unless we impose strong growth conditions on $w$ at zero. We therefore express the solution of \eqref{y1} by means of the modified function

\begin{align}\label{phi}
\phi(\delta) = \int_0^\delta \frac{y_1}{s}\:ds \hspace{90pt} \\
 = c \delta^{-2\beta}+\frac{A\delta}{2\beta +1} - \delta^{-2\beta}\int_\delta^{\delta_0} Hs^{2\beta-1}\:ds.
\end{align}

In this case the leading term of $\phi$ is of order $O(\delta)$ provided $A\neq 0$. If $A=0$ it is of higher order. We consider the operator $\Theta(w)$ in the Banach space
$X:= \{ w\in C([0,\delta_0]): |w|\leq M\delta^2\:, \,\,\delta \in [0,\delta_0]\}$, where $M$ and $\delta_0\le \frac{R-r_0}2$ are positive constants which will be determined later, and $\|w\|:= \sup\{|\frac{w(\delta)}{\delta^2}|\:, \,\,\delta \in (0,\delta_0]\}$.
The estimates are similar to the ones in the first case except that
\begin{align*}
|H(\delta,w)-H(\delta,\tilde w)| \leq c_4\delta^2\|w-\tilde w\|,\\
|\phi(\delta)-\tilde \phi(\delta)|\leq c_5\delta^2\|w-\tilde w\|.
\end{align*}
If $\beta \not= -1$, as before this leads to $\|\theta(w)-\Theta(\tilde w)\|\leq c_6 \delta_0\|w-\tilde w\|$. Notice that if $\beta \not= -1$ the expression
$\delta^{-2\beta}\int_\delta^{\delta_0} Hs^{2\beta-1}\:ds$ is of order $O(\delta^2)$.
For the next claim that $\Theta: X\to X$ we observe that for $w\in X$
$$
|T(w)|\leq c \delta^{-2\beta} +\frac{|A|}{|2\beta +1|} \delta +c_2M \delta^2.
$$
Then
$$\|\Theta(w)\| \le c \delta_0^{-2\beta -1} + \frac{|A|}{2|2\beta +1|}  +c_2 M \delta_0 <M\,, $$
for $M > \frac{|A|}{2|2\beta +1|}$ and $\delta_0$ sufficiently small.
From here we conclude that $\Theta:X\to X$ is a contraction and has a unique fixed point.\par
If $\beta = -1$, then $\delta^{-2\beta}\int_\delta^{\delta_0} Hs^{2\beta-1}\:ds$ is of order $O(\delta^2|\ln \delta|)$. By requiring that $\delta_0|\ln \delta_0|$  is sufficiently small we obtain that $\Theta$ is a  contraction in $X$.
\par
Notice that the dependence of the constant $C$ from $c$ is not explicit in this case.
\medskip

\hfill $\square$
\medskip

\begin{remark} The constant $A$ vanishes if $\beta=-\frac{N-1}{N+1}$ or if $N=1$. If both $c$ and $A$ vanish higher order terms come into play. The discussion is straightforward and will be omitted.
\end{remark}
\section{Global  solutions}
\subsection{Ball}
\begin{theorem}\label{thm:ball} Assume $\mu<1/4,\, \mu \not=0$. For $\Omega =B_R$ we have

(i) For any given $u(0)> 0$ problem \eqref{original} possesses in the ball a unique positive radial solution. At the boundary it behaves like $c\delta^{\beta_-}$, for some $c>0$.  The solutions are monotone increasing with respect to $u(0)$.

(ii) For any $0< R_0<R$ there exists a  nonnegative radial solution in the ball with a dead core in $B_{R_0}$. At the boundary it behaves like $c\delta^{\beta_-}$ for some $c>0$.

(iii) There exists a solution of the form $u(r)= r^{\frac{2}{1-p}}(c'' +w(r))$ with $c''=\big(-\mu^*+\frac{2(N-1)}{1-p} \big)^{\frac1{p-1}}$ and $w(0)=0$. At the boundary it behaves like $c\delta^{\beta_-}$, for some $c>0$.

\end{theorem}
{\it Proof.} From Section \ref{locsol} we know that problem \eqref{Er} with the initial conditions $u(0)=u_0>0$ and $u'(0)=0$ has a unique local solution which can be continued until it vanishes or it blows up. Since $p<1$ blow up cannot occur for $r<R$. If $\mu<1/4$ then by the comparison principle stated in the introduction it cannot vanish before $r=R$.
By the results of the previous section it behaves at the boundary like $c\delta^{\beta_-}$  with $c>0$ or it is bounded from above by $c\delta^{\beta_+}$. The second case is impossible in view of the comparison principle. Consequently $u\sim c\delta^{\beta_-}$ at the boundary. Solutions are monotone increasing with respect to $u(0)$ since they cannot intersect for $r \in (0, R)$ by the comparison principle.

These solutions are positive in the whole ball. All other solutions have a dead core. In fact if we choose $R_0>0$,  set $u=0$ in $[0,R_0]$ and continue it with the solution constructed in Section \ref{locsol}, by the same arguments as before we obtain a solution which exists in the whole ball and behaves at the boundary like $c\delta^{\beta_-}$. Notice that a solution for which $u(R_0)=0$ and $u'(R_0)=0$ is necessarily zero in $(0,R_0)$.

The third assertion follows from Lemma \ref{deadcore} (iii). \hfill $\square$
\medskip

If $\mu<0$ the solutions are monotone increasing and blow up at the boundary. This is not the case if $\mu>0$. 

Notice that the solution with a dead core at the boundary has a singularity at the origin.
\subsection{Annulus}
The structure of the positive {\it radial solutions}  in an annulus is described in
\begin{theorem}\label{thm:annulus}
For $\Omega = \mathcal{A}(r_0,R)$  and $\mu<C_H(\mathcal{A}(r_0,R)),\, \mu \not=0$, we have

(i) For any given $r_0<R_0 <R$ there exists a unique solution positive in $(R_0,R)$, with a dead core in $[r_0,R_0]$. At the outer boundary it behaves like $k (R-r)^{\beta_-}$, for some $k>0$.  Vice versa

For any given $r_0<R_0 <R$ there exists a unique solution positive in $(r_0,R_0)$,  with a dead core in $[R_0,R]$. At the inner boundary it behaves like $k(r-r_0)^{\beta_-}$, for some $k>0$.

(ii) The sum of two solutions as in (ii), having a disjoint support, is a solution with a dead core interval (eventually reduced to a point) and positive near the inner and outer boundary.

(iii) If $\mu <0$, for any given $r_0<R_0 <R$ and $u(R_0) =u_0 >0$ there exists a unique positive solution. At the outer and inner boundaries it behaves like $k_1 (R-r)^{\beta_-}$, respectively $k_2(r-r_0)^{\beta_-}$, for some $k_1, k_2 >0$.

If moreover  $\mu> -\frac{2(p+1)}{(1-p)^2}$ we have

(iv) For any given $c >0$ there exists a unique positive solution such that $u(r)/(r-r_0)^{\beta_+}\to c$ as $r\to r_0$. At the outer boundary it behaves like $k(R-r)^{\beta_-}$, for some $k>0$. Vice versa we have

for any given $c>0$ there exists a unique positive solution such that $u(r)/(R-r)^{\beta_+}\to c$ as $r\to R$. At the inner boundary it behaves like $k(r-r_0)^{\beta_-}$, for some $k>0$.

 (v)  There exists a unique solution such that
 \newline $u(r)/(r-r_0)^{\frac{2}{1-p}}\to c'$ as $r\to r_0$. At the outer boundary it behaves like
 \newline $k(R-r)^{\beta_-}$, for some $k>0$. Here
$ c':=\left( \frac{2(1+p)}{(p-1)^2} +\mu\right)^{1/(p-1)}$.

 The same holds if we interchange  the role of the inner and outer boundary.
\end{theorem}
{\it Proof.}  In order to prove the first statement let $\tilde u$ be the solution with a dead core in one point $r=R_0$ constructed in Lemma \ref{deadcore}. The same arguments that we used in the proof of Theorem \ref{thm:ball} give that $\tilde{u}$ can be continued to the right and the left until it reaches the inner and outer boundary. There it behaves like $k\delta^{\beta_-}$, where $\delta$ denotes the distance from the boundary and $k>0$. If in $[r_0,R_0]$ (respectively in $[R_0,R]$) we replace it with $\tilde{u} \equiv 0$, we get (i).

(ii) is a simple remark.

(iii) As already remarked, problem \eqref{Er}, \eqref{ic2} has a local solution. Moreover for $u_1=0$ and $\mu <0$ this solution increases in $[R_0, R)$ and decreases in $(r_0, R_0]$, hence it is positive and cannot go to $0$ at the boundary. By Corollary \ref{Corb-} the solution behaves as $k\delta^{\beta_-}$ at the inner and outer boundary. By Lemma \ref{le:b+} we have $k>0$.

(iv) we start with the local solution which behaves at the inner or outer boundary like
$c\delta^{\beta_+}$ (see Lemma \ref{localregular}). It can be continued till the outer or inner boundary. Then we argue as in Theorem \ref{thm:ball}.

(v) is proved exactly on the same line (see Lemma \ref{deadcore}, (ii)).\hfill $\square$
\
\subsection{General domains}
In this section we shall construct solutions of \eqref{original} in arbitrary not necessarily simply connected domains.  More precisely we shall prove the following theorem.
\begin{theorem}\label{thGD} Let $\mu<\frac14,\, \mu\not=0,$ and $\Omega$ be a bounded domain with $C^k$ ($k\geq 2$) boundary. Then the following statements hold for the solutions of \eqref{original}:
\begin{description}
\item[(i)]  for suitable $0< c_0 < c_1$ there exists a solution $u$ such that
$0<c_0\leq u(x)/\delta^{\beta_-}(x) \leq c_1$ in a neighborhood of $\partial\Omega$.

\item[(ii)] If $c_0$ and $c_1$ are sufficiently small this solution has a dead core in the interior of $\Omega$.

\item[(iii)] If $\partial\Omega$ is not connected, then for any non empty, closed, disjoint sets $\Gamma_1,\,\Gamma_2$, such that $\Gamma_1\cup\Gamma_2 = \partial\Omega$, and for suitable sufficiently small $0< c_0 < c_1$, there exists a solution $u$ positive in a neighborhood of $\Gamma_1$ where it behaves as in $(i)$
and such that $u(x) \equiv 0$, in a neighborhood of $\Gamma_2$.
\end{description}
\end{theorem}


For the proof of the theorem we need some properties of the distance function $\delta(x)$ where $x$ is an arbitrary point in $\Omega$.  Denote by $\Omega_\rho$ the parallel set $\{x\in \Omega: \delta(x)<\rho\}$. If $\Omega$ is of class $C^k$, $k\geq 2$, then $\delta$ is in $C^k(\Omega_{\rho_0})$ for $\rho_0>0$ sufficiently small. Denote by $\sigma(x)$ the nearest point to $x$ on $\partial \Omega$. Let $K_i(\sigma(x))$, $i=1,..,N-1$  be the principal curvatures and  $H(\sigma(x))=\sum_{i=1}^{N-1}\frac{K_i}{N-1}$  be the mean curvature. Then for any $x \in \Omega_{\rho_0}$
\begin{equation}\label{Deltadelta}\begin{split}
|\nabla \delta (x)| =1\,,\hspace{90pt}\\
 -\frac{N-1}{\rho_0 - \delta(x)}\le \Delta \delta(x) = - \sum_{i=1}^{N-1}\frac{K_i}{1-K_i\delta(x)} \le \frac{N-1}{\rho_0 + \delta(x)}\,.
 \end{split}
\end{equation}
\bigskip

{\it Proof of Theorem \ref{thGD}}.

$(i)$ For the proof  of the first assertion we shall distinguish between two cases.

{\bf (A)}\, $\mu\in (0,1/4)$.
\medskip

\noindent For $0 < s \le \rho < \frac{\rho_0}2$, $\epsilon >0$, let $\phi(s):= Ms^{\beta_-}(\rho^{\epsilon} - s^{\epsilon})$.  Then
\begin{align*}
 \phi'(s)= \beta_- Ms^{\beta_--1}(\rho^{\epsilon} - \frac{\beta_- + \epsilon}{\beta_-} s^{\epsilon}),\\
  \phi''(s):= \beta_-(\beta_--1) \frac{\phi(s)}{s^2} - M \epsilon (2 \beta_- + \epsilon -1)s^{\beta_- +\epsilon-2} .
  \end{align*}
The function $\tilde u(x):= \phi(\delta(x))$ is well defined  for $x \in \Omega_{\rho_0}$,  and it satisfies (by \eqref{Deltadelta}):
\begin{equation}\label{phididelta}\begin{array}{c}
\Delta \tilde u(x) = \phi''(\delta)|\nabla \delta|^2+ \phi'(\delta) \Delta \delta  \\
= -\mu \frac{\phi(\delta)}{\delta^2}  - M \epsilon (2 \beta_- + \epsilon - 1)\delta^{\beta_- +\epsilon-2} +
\phi'(\delta) \Delta \delta\,.
\end{array}
\end{equation}
By a suitable choice of $\epsilon$ we can construct {\it local upper and lower solutions}. In fact:
\begin{itemize}
  \item[(a)] if $0< 1-2\beta_- <\epsilon < 1$, then there exists $\rho<\frac{\rho_0}2$ sufficiently small such that, for any $M>0$, $\tilde u$ is an upper solution in $\Omega_{\rho}$.
  \item[(b)] For any given $0<\epsilon:=\underline{\epsilon}<1-2\beta_-<1$ and $M>0$ there exists $\rho<\frac{\rho_0}2$ sufficiently small such that $\tilde u$ is a lower solution in $\Omega_{\rho}$.
\end{itemize}
The first assertion (a) follows
from the estimate
$$
 |\phi'(s)| \le \beta_-Ms^{\beta_--1}\rho^{\epsilon}\max\{1,\frac{\epsilon}{\beta_-}\}\le M K s^{\beta_--1}\,
$$
for some constant $K$ independent of $s$, and
$$
|\Delta \delta(x)| \leq K_1 \tx{in} \Omega_\rho,
$$
where $K_1$ depends only on $\rho_0$. Inserting these estimates into \eqref{phididelta} we get
\begin{align*}
\Delta \tilde{u} + \mu \frac{\tilde{u}}{\delta^2} - \tilde{u}^p \le - \epsilon(2 \beta_- + \epsilon-1)M\delta^{\beta_- +\epsilon-2} + MK K_1 \delta^{\beta_--1}  \\
=  - M\delta^{\beta_- +\epsilon-2}\left[\epsilon (\epsilon-(1-2 \beta_-)) - KK_1 \delta^{1 - \epsilon}\right].
\end{align*}
For small $\delta$ the right-hand side is negative. This proves the first assertion.
\medskip

The second assertion $(b)$ follows from
\begin{align*}
\Delta \tilde{u} + \mu \frac{\tilde{u}}{\delta^2} - \tilde{u}^p \ge  M [\underline{\epsilon} (1-2 \beta_- - \underline{\epsilon})\delta^{\beta_- +\underline{\epsilon}-2} -K K_1 \delta^{\beta_--1} - K_2M^{p-1}\delta^{p\beta_-}] \\
=   M \delta^{\beta_- +\underline{\epsilon}-2}[\underline{\epsilon} (1-2 \beta_- -\underline{\epsilon}) - K K_1 \delta^{1 - \underline{\epsilon}} - K_2M^{p-1} \delta^{2-(1-p)\beta_-- \underline{\epsilon}}] >0\,,
\end{align*}
where $K_2$ depends only on $\rho_0$ and $\epsilon$. Since $2-(1-p)\beta_- -\epsilon >0$ the right-hand side is positive for small $\rho$. This completes the proof of $(b)$.
\bigskip

Next we want to extend the local upper and lower solutions constructed above to the whole domain. Let $\rho \in (0, \frac{\rho_0}2]$ be such that
$$
\bar{u}= M\delta^{\beta_-}(\rho^\epsilon-\delta^\epsilon)
$$
is an upper  solution in $\Omega_{\rho}$.

Observe that $\bar{u}$ attains its maximum $\bar{u}_M$ at $\{ x \in \Omega \, \,: \delta(x) = \bar{\delta}:= (\frac{\beta_-}{\beta_- + {\epsilon}})^\frac1{{\epsilon}}\rho\}$.\par
We choose $M$ so small that the following inequality holds:
\begin{equation}\label{chooseM}
\mu \frac{\overline{u}_M}{\bar{\delta}^2} - \overline{u}_M^p = \overline{u}_M^p [\mu \frac{\overline{u}_M^{1-p}}{\bar{\delta}^2} - 1]<0\,.
\end{equation}
Then the constant $\overline{u}_M$ is an upper solution of \eqref{original} in $\Omega \setminus \Omega_{\bar{\delta}}$ and we obtain the following (weak) global upper solution
\begin{equation}\label{baru}
\bar{U}(x): = \left\{\begin{array}{c}
            \overline{u}(x)\,,\,\,x \in \Omega_{\bar{\delta}}\,,
             \\
            \overline{u}_M \,,\,\,x \in \Omega \setminus \Omega_{\bar{\delta}}\,.
          \end{array}\right.
\end{equation}

For the same $M$ let $\underline{\rho} \in (0, \bar{\delta})$ be such that  $\underline{u}=M \delta^{\beta_-}(\underline{\rho}^{\underline{\epsilon}}-\delta^{\underline{\epsilon}})$ is a lower solution in $\Omega_{\underline{\rho}}$, such that $\bar{u}>\underline{u}$ in $\Omega_{\underline{\rho}}\subset \Omega_{\bar{\delta}}$. The function
$$
\underline{U}(x)=
\begin{cases}
\underline{u} &\tx{in} \Omega_\rho,\\
0 &\tx{otherwise}
\end{cases}
$$
is a global lower solution.

Hence there exist an upper and a lower solution $\underline{U}(x)\leq \overline{U}(x)$ in $\Omega$. The method of upper and lower solutions can be generalized to our case cf. \cite{BaMoRe08} (Lemma 4.12) and guarantees the existence of a non trivial positive solution $\underline{U}\leq u\leq \overline{U}$.
\bigskip

{\bf (B)}\, $\mu<0$.
\medskip

\noindent We start with the construction of an upper solution. Let $\bar{\sigma}$ be a positive number smaller than $\rho_0$ and for any given $M>0$ let $\eta=\eta(r)$ be the solution of
\begin{equation*}
\left\{\begin{array}{c}
\eta'' +\frac{(N-1)}{r}\eta' + \frac{\mu}{(\rho_0-r)^2}\eta = \eta^p\,,\,\, r \in (\rho_0 - \bar{\sigma}, \rho_0),\\
\eta(\rho_0-\bar{\sigma})=M\,,\,\, \eta'(\rho_0-\bar{\sigma})=0\,.  \hspace{35pt}                                                                                    \end{array}\right.
\end{equation*}
Since $\mu<0$ the function $\eta(r)$ is increasing to the right and can therefore  be extended as a positive solution in $(\rho_0-\bar{\sigma},\rho_0)$. Then
\begin{equation}\label{bryUn}
\lim_{r \to \rho_0} \frac{\eta(r)}{(\rho_0-r)^{\beta_-}} = C_M >0 \,.
\end{equation}
Indeed by Lemma \ref{le:b+} $(ii)$, if $C_M=0$ then $\eta(r) \to 0$ as $r \to \rho_0$ which contradicts the increasing behavior of $\eta$.\par
We can easily verify that the following function
\begin{equation}\label{barun}
\bar{u}(x): = \left\{\begin{array}{c}
            \eta(\rho_0 -\delta(x))\,,\,\,x \in \Omega_{\bar{\sigma}}\,,
             \\
            M \,,\,\,x \in \Omega \setminus \Omega_{\bar{\sigma}}\,.
          \end{array}\right.
\end{equation}
is a (weak) upper solution of \eqref{original}. Indeed since $\mu <0$, any constant is an upper solution. Since  $\bar{u} \in C^1(\Omega)$,
we only have to verify that it is a classical upper solution for any $x \in \Omega_{\bar{\sigma}}$.
Indeed remark that $\eta'$ is positive, hence for any $x \in \Omega_{\bar{\sigma}}$ we have by \eqref{Deltadelta} and \eqref{phididelta}
\begin{equation}\label{DeltaUn}\begin{split}
 \Delta \bar{u}(x) = \eta''(\rho_0 -\delta(x))|\nabla \delta|^2 - \eta'(\rho_0 -\delta(x))\Delta\delta(x) \\ \le  \eta''(\rho_0 -\delta(x)) + \frac{(N-1)}{\rho_0 -\delta(x)}\eta'(\rho_0 -\delta(x)) \\ = - \frac{\mu}{\delta(x)^2}\eta + \eta^p  = - \frac{\mu}{\delta(x)^2}\bar{u} + \bar{u}^p \,.
 \end{split}
\end{equation}
\medskip

In order to construct a lower solution take $\underline{\sigma}\in (0,\rho_0)$   and let $z=z(d)$ be the non trivial (dead core) solution of
\begin{equation}\label{zeta}
\left\{\begin{array}{c}
z'' +\frac{(N-1)}{\rho_0 + d}z' + \frac{\mu}{d^2}z = z^p\,,\,\,d \in (0,\underline{\sigma}) \\
z(\underline{\sigma})=0\,,\,\, z'(\underline{\sigma})=0\,.  \hspace{35pt},\:                                                                                    \end{array}\right.
\end{equation}
such that $z(d)>0\,, \,\,d\in (0, \underline{\sigma})$. We extend it by $0$ for $d \ge \underline{\sigma}$, set $w(r):=z(r-\rho_0)$ and observe that it is a radial solution of \eqref{original} in the annulus $\mathcal{A}(\rho_0, R)$, for any $R>\rho_0+2\underline{\sigma}$.  In addition
\begin{equation}\label{bryUnl}
\lim_{d \to 0} \frac{z(d)}{d^{\beta_-}} = C_{\underline{\sigma}} >0 \,.
\end{equation}
We can easily verify that the following function
\begin{equation}\label{underlineun}
\underline{u}(x): = \left\{\begin{array}{c}
            z(\delta(x))\,,\,\,x \in \Omega_{\underline{\sigma}}\,,
             \\
           0 \,,\,\,x \in \Omega \setminus \Omega_{\underline{\sigma}}\,.
          \end{array}\right.
\end{equation}
is a (weak) lower solution of \eqref{original}. Indeed  $\underline{u} \in C^1(\Omega)$ and it satisfies \eqref{original} in the classical sense in the interior of the region where it vanishes. Hence
we only have to verify that it is a classical lower solution for any $x \in \Omega_{\sigma}$.
Indeed remark that $z'$ is negative, hence for any $x \in \Omega_{\underline{\sigma}}$, by \eqref{phididelta} we have
\begin{equation}\label{DeltaUnl}\begin{split}
 \Delta \underline{u}(x) = z''(\delta(x))+ z'(\delta(x))\Delta\delta(x) \\ \ge  z''(\delta(x)) + \frac{(N-1)}{\rho_0 +\delta(x)}z'(\delta(x)) \\ = - \frac{\mu}{\delta(x)^2}z + z^p  = - \frac{\mu}{\delta(x)^2}\underline{u} + \underline{u}^p \,.
 \end{split}
\end{equation}
It is not difficult to see that by choosing $\underline{\sigma}$ sufficiently small we can achieve that $\underline{u}\leq \overline{u}$. Hence the proof is complete.\bigskip


\noindent
$(ii)$ We distinguish between two cases as in $(i)$.
\bigskip

{\bf (A)}\, $\mu \in (0, \frac14)$.
\medskip

We construct an upper (weak) solution with dead core. Let $\bar{U}(x)$ be the upper solution in $(i)$, {\bf (A)}, and $\bar{\delta}\in (0, \rho_0)$ the constant used in its definition. Take $\rho \in
(\bar{\delta}, \rho_0)$.
By Lemma \ref{deadcore} there exists $\eta=\eta(r)$ solution of
\begin{equation}\label{eta}
\left\{\begin{array}{c}
\eta'' +\frac{(N-1)}{r}\eta' + \frac{\mu}{(\rho_0-r)^2}\eta = \eta^p\,,\,\, r \in (\rho_0-\rho, \rho_0-\bar{\delta}),\\
\eta(\rho_0-\rho)= 0\,,\,\, \eta'(\rho_0-\rho)=0\,.  \hspace{35pt}                                                                                    \end{array}\right.
\end{equation}
For $\eta$ sufficiently close to $\eta=0$ and $ r \in (\rho_0-\rho, \rho_0-\bar{\delta})$, the quantity $\eta^p(1-\frac{\mu}{(\rho_0-r)^2}\eta^{1-p})$ is positive, hence $\eta(r)$ is increasing to the right in a small interval $(\rho_0 - \rho, \rho_0 - \tilde{\rho})$, for some $\tilde{\rho}\in (\bar{\delta}, \rho)$.
 As in \eqref{DeltaUn} we obtain that $\eta(\rho_0-\delta(x))$ is a local upper solution of \eqref{original} in $\Omega_{\rho} \setminus \Omega_{\tilde{\rho}}$.\par
If $\bar{U}(x)\le \eta(\rho_0-\tilde{\rho})$ for $\delta(x)= \tilde{\rho}$, there exists an eventually larger $\tilde{\rho}$, such that $\bar{U}(x) = \eta(\rho_0-\tilde{\rho})$ for $\delta(x)= \tilde{\rho}$.\par
The following function is a (weak) global upper solution  with dead core
\begin{equation*}
\tilde{U}(x): = \left\{\begin{array}{c}
            \overline{U}(x)\,,\,\,x \in \Omega_{\tilde{\rho}}\,,
             \\
             \eta(\rho_0-\delta(x))\,,\,\,x \in  \Omega_{\rho} \setminus\Omega_{\tilde{\rho}}\,,
             \\
            0 \,,\,\,x \in \Omega \setminus \Omega_{\rho}\,.
          \end{array}\right.
\end{equation*}
If $\bar{U}(x)> \eta(\rho_0-\tilde{\rho})$, we remark that for any $m \in (0,1)$, $m \bar{U}(x)$ is an upper solution,
hence we can choose $m$ such that $m\bar{U}(x)= \eta(\rho_0-\tilde{\rho})$ for $\delta(x)= \tilde{\rho}$ and conclude as above.\par
Concerning the lower solution we proceed as in $(i)$,{\bf (A)}.
 The conclusion follows as in case (i).
\bigskip

{\bf (B)}\, $\mu<0$.
\medskip

We construct the upper solution as in $(i)$, {\bf (B)}\, solving \eqref{eta} for $M=0$. A positive solution exists by Lemma \ref{deadcore}. No other changes are needed in the proof.

\bigskip

\noindent
$(iii)$ For any $\rho \in (0, \rho_0)$ we define the following subsets of $\Omega$, say $(\Gamma_1)_{\rho}:=\{x \in \Omega\, :\, d(x,\Gamma_1)<\rho\} $ and $(\Gamma_2)_{\rho}:=\{x \in \Omega\, :\, d(x,\Gamma_2)<\rho\} $.  $\rho_0$ is such that  $(\Gamma_1)_{\rho_0}$ and $(\Gamma_2)_{\rho_0}$ are disjoint sets and
 $\Omega_{\rho_0} =(\Gamma_1)_{\rho_0} \cup (\Gamma_2)_{\rho_0}$.\par
  In (ii) we constructed a solution $u$ which vanishes in $\Omega \setminus \Omega_{\rho_0}$, hence the function
  \begin{equation*}
\tilde{u}(x): = \left\{\begin{array}{c}
            u(x)\,,\,\,x \in (\Gamma_1)_{\rho_0}\,,
             \\
            0 \,,\,\,x \in \Omega \setminus (\Gamma_1)_{\rho_0}\,,
          \end{array}\right.
\end{equation*}
is a solution and it has the behavior required in (iii).
\hfill $\square$
\begin{remark}
\begin{enumerate}
\item In the proof of $(i)$ we have constructed upper solutions and smaller nontrivial lower solutions which vanish in an interior set. The solutions we constructed lie between the upper and lower solutions, hence they might be strictly positive or  have a dead core in some subsets of $\Omega$.

\item Under the hypotheses in $(iii)$, the number of solutions with a different qualitative behavior depends on the possible choices of the sets $\Gamma_1$ and $\Gamma_2$, hence on the number of connected components of $\partial \Omega$. In particular the role of $\Gamma_1$ and $\Gamma_2$ can be exchanged.

\item We only required $\mu <\frac14$ and not $\mu <C_H(\Omega)$. This weaker requirement is due to the fact that the Hardy constant of a thin set is $\frac14$ and we mainly work in a thin neighborhood of the boundary.\end{enumerate}
\end{remark}


\end{document}